\documentclass[11pt]{amsart}
\usepackage{epsf} 
\usepackage{amsmath, amsthm, amssymb, amsfonts, verbatim, color}
\usepackage[pdftex]{graphicx}
\usepackage[bookmarks]{hyperref}

\renewcommand{\baselinestretch}{1}
\reversemarginpar

 \newcommand{\norm}[1]{\left\Vert#1\right\Vert}
\newcommand{\ip}[1]{\langle{#1}\rangle}

\newcommand{\abs}[1]{\left\vert#1\right\vert}

\def\be{\begin{equation}}
\def\ee{\end{equation}}
\def\qand{\quad\mbox{and}\quad}

\def\Ha{{\mathcal H^{s',\theta_{1}}}}
\def\H{{\mathcal H^{s,\theta}}}
\def\h{{H^{s,\theta_0}}}

\newcommand{\dx} {\; \mathrm{d} x}
\newcommand{\dd} {\mathrm{d}}

\def\a{\alpha}
\def\b{\beta}
\def\eps{\epsilon}

\DeclareMathOperator{\diag}{diag}

\DeclareMathOperator{\dv}{div}
\DeclareMathOperator{\BOX}{\square}

\def\R{{\mathbb R}}
\def\C{{\mathbb C}}
 
\def\bs{\begin{split}}

\theoremstyle{plain}

\def\lm{\begin{lem}}
\def\ml{\end{lem}}
\newtheorem{thm}{Theorem}[section]
\newtheorem{cor}[thm]{Corollary}
\newtheorem{lem}[thm]{Lemma}

\theoremstyle{definition}

\theoremstyle{remark}
\newtheorem{remark}{Remark}[section]

\numberwithin{equation}{section}

\def\be{\begin{equation}}
\def\ee{\end{equation}}
\def\qand{\quad\mbox{and}\quad}

\def\R{\mathbb R}

\def\pru{\begin{proof}}
\def\urp{\end{proof}}

\theoremstyle{definition}

\newtheorem*{theorem*}{Theorem}
\newtheorem*{corollary*}{Corollary}
\newtheorem*{acknowledgment}{Acknowledgments}

\begin{document}

\date{\today}

\title[2D MKG in Coulomb]{Low regularity well-posedness for the 2D Maxwell-Klein-Gordon equation in the Coulomb gauge}
\author[Czubak]{Magdalena Czubak}
\email{czubak@math.binghamton.edu}
\author[Pikula]{Nina Pikula}
\email{npikula1@binghamton.edu}
\address{Department of Mathematical Sciences, Binghamton University (SUNY)}

\begin{abstract}
We consider the Maxwell-Klein-Gordon equation in 2D in the Coulomb gauge.  We establish local well-posedness for $s=\frac 14+\epsilon$ for data for the spatial part of the gauge potentials and for $s=\frac 58+\epsilon$ for the solution $\phi$ of the gauged Klein-Gordon equation.   The main tool for handling the wave equations is the product estimate established by D'Ancona, Foschi, and Selberg.   Due to low regularity, we are unable to use the conventional approaches to handle the elliptic variable $A_0$, so we provide a new approach.
\end{abstract}

\maketitle
%%%%%%%%%%%%%%%%%%%%%%%%%%%%%%%%%%%%%%%%%%%%%%%%%
%%%%%%%%%%%%%%%%%%%%  Introduction
%%%%%%%%%%%%%%%%%%%%%%%%%%%%%%%%%%%%%%%%%%%%%%%%%

\section{Introduction}\label{intro}

We study local well-posedness (LWP) of the Cauchy problem for the 2D Maxwell-Klein-Gordon equation (MKG) in the Coulomb gauge.    Well-posedness for MKG in 2D has been so far only considered in the Lorenz and temporal gauges.  Moncrief  \cite{MoncriefMKG} showed global well-posedness in the Lorenz gauge  for data in $H^2$.  Recently Pecher \cite{Pecher} studied LWP for data with $s=\frac 14+\epsilon$ for the gauge potentials and $s=\frac 34+\epsilon$ for $\phi$, the solution of the gauged Klein-Gordon equation.  In the temporal gauge, there is work by Schwarz \cite{Schwarz} for $s\geq 2$ and with $\abs{\phi}\rightarrow 1$ at infinity.

Based on the previous works for wave equations in 2D the common expectation could be that MKG in the Coulomb gauge should be well-posed for $s>\frac 14$ (we explain this below).  Moreover, it might seem that this is obvious and that it should simply follow from well-known estimates.  However, at this low level of regularity even solving the elliptic equation comes with obstacles.  As a result, low regularity well-posedness for MKG in 2D becomes more interesting than initially expected.
 
In the Coulomb gauge, MKG is a system of wave equations for the complex field $\phi$ and the spatial part of the connection $A$ coupled to elliptic equations for the temporal parts, $A_{0}$ and $\partial_{t}A_{0}$.  The nonlinearities involve null forms and other bilinear and trilinear terms.   

The null condition was introduced by Klainerman in \cite{Klainerman83}, and it was first used to lower regularity assumptions on initial data in \cite{KlainermanMachedon93}.   The null form in MKG in the Coulomb gauge was originally uncovered in   \cite{KM}.
% and then written for all dimensions in \cite{SelbergMKG}.  

The null form appearing in MKG is
\[
Q_{\alpha\beta}(u,v)=\partial_\alpha u \partial_\beta v - \partial_\beta u \partial_\alpha v.
\]
In 3D, almost optimal LWP for initial data in $H^s\times H^{s-1}, s>\frac 32$, for wave equations with this particular null form was shown   in \cite{KM3}.   For example, for systems that can be written as
\[
\Box u=Q_{\alpha\beta}(u,v)=\Box v, \quad \alpha,\beta \in \{0,\dots, 3\}, \ \alpha\neq \beta.
\]
In 2D, the situation is not as optimal.    Note that by scaling invariance, $s=\frac n2$ is the critical exponent for the system on $\R^{n+1}$.  By examining the first iterate Zhou \cite{Zhou} showed that $s>\frac 54$ is as close as one can get using iteration methods (but see \cite{GN}). 

Now, if we do not consider the elliptic equation, cubic terms or bilinear terms involving the elliptic variables, MKG in the Coulomb gauge can be schematically written as
\be\label{model}
\begin{split}
\Box u&=D^{-1}Q_{\alpha\beta}(u,v),\\
\Box v&=Q_{\alpha\beta}(D^{-1}u,v).
\end{split}
\ee
The presence of $D^{-1}=(-\Delta)^{-\frac 12}$ changes the scaling transformation and shifts the critical exponent to $\frac n2-1$.

Machedon and Sterbenz  \cite{MachedonSterbenz} established almost optimal LWP for MKG in 3D for $s>\frac 12$.  In addition, they showed that the system \eqref{model} will be \emph{ill-posed} below $\frac 34$ if one only considers the above model equations.

Now, because of \cite{Zhou} and $D^{-1}$ heuristically one might expect MKG to be locally well-posed for $s>\frac 14$.  However, this same heuristic raises an expectation of \eqref{model} being well-posed for $s>\frac 12$ in 3D, but again \cite{MachedonSterbenz} showed $s>\frac 34$ is needed.

In this paper we show that $s>\frac 12$ is needed if we only use the wave analog of $X^{s,b}$ spaces (defined in Section \ref{spaces} below) and assume $A$ and $\phi$ have the same regularity.   However, we also show that we can let $s'=\frac 14+\epsilon$ for $A_j$ if $s=\frac 58+\epsilon$ for $\phi$ (see Theorem \ref{mainthm} for a precise statement).

We note that one of the observations that allowed \cite{MachedonSterbenz} to lower the regularity was the recognition of the cancellations between the null form and the elliptic term in the wave equation for $\phi$.  Here the need for $s>\frac 12$ comes already in the equation for $A_j$.

The main technical tool for handling the wave equation estimates is the convenient atlas of product estimates established by D'Ancona, Foschi, and Selberg \cite{DAFS} (see Theorem \ref{atlas} below).

Finally, controlling the elliptic estimates for $A_0(t)$ in 2D when $s<1$ causes difficulties when one attempts the standard methods.  We provide an alternative approach to resolve this.

\subsection{MKG system and the statement of the results}
MKG is a system of Euler-Lagrange equations of the following action functional
\be\label{lag}
\mathcal L (A, \phi)=\int_{\R^{2+1}}\frac 12 D_{\a}\phi D^{\a}\phi+\frac 14F_{\a\b}F^{\a\b}\dx \dd t.
\ee
Here
\[
\phi: \R^{2+1}\to \C,\quad
A\mbox{ is  a $1$-form with components }\quad A_\a: \R^{2+1}\rightarrow \R
, \ \ \alpha\in \{0, 1, 2\},
\]
$D_\alpha$ denotes the covariant derivative
$$
D_{\alpha}\phi := (\partial_{\alpha}+iA_{\alpha})\phi,
$$
and $F := dA$, so that
\[
F_{\a\b} = \partial_\a A_\b - \partial_\b A_\a, \qquad\a,\b\in \{0, 1, 2\}.
\]
One may regard $A$ as a $U(1)$ connection and $F$ as the associated curvature.

In \eqref{lag} we  sum over repeated upper and lower indices, and
we raise and lower indices with the Minkowski metric
$(\eta^{\a\b}) = (\eta_{\a\b}) = \mbox{diag}(-1,1,1)$
so that
\[
D^{\alpha}\phi = \eta^{\a\b}D_{\b}\phi,
\quad\quad
F^{\a\b} = 
\eta^{\a\gamma} \eta^{\b \delta} F_{\gamma\delta}.
\]
The Euler-Lagrange equations associated with the action functional \eqref{lag} are
\begin{subequations}
\begin{align}
D_{\a}D^{\a}\phi&=0,\label{mkg1}\\
\partial_{\a}F^{\a\b}&=J^{\b},\quad \b\in \{0,1,2\}\label{mkg2},
\end{align}
\end{subequations}
where $J^{\b}$ is the current given by
\[
J^{\b}=-\Im (\phi\overline{ D^{\b}\phi}),
\]
and $\Im(z)$ denotes the imaginary part of the complex number $z$.
 
The action functional \eqref{lag}  is invariant under the action of the $U(1)$ group, so for any sufficiently regular $f:\R^{2+1}\rightarrow \R$ we have
\[
\mathcal L(\phi,A)=\mathcal L( e^{if} \phi, A- df).
\]
In the Coulomb gauge, $\partial^{j}A_{j}=0$, and \eqref{mkg1}-\eqref{mkg2} become
\begin{subequations}
\begin{align}
\Delta A_{0}&=-\Im (\phi\overline{ \partial_{t}\phi})+\abs{\phi}^{2}A_{0},\label{mkga0}\\
\BOX A_{j}&=-\Im (\phi\overline{ \partial_{j}\phi})+\abs{\phi}^{2}A_{j}-\partial_{j}\partial_{t}A_{0},\label{mkgaj}\\
\BOX \phi&=-2iA^{j}\partial_{j}\phi +2iA_{0}\partial_{t}\phi+i(\partial_{t}A_{0})\phi+A^{\a}A_{\a}\phi,\label{mkgphi}\\
\partial^{j}A_{j}&=0\label{ccond}.
\end{align}
\end{subequations}
As is now well-known, the equations \eqref{mkga0}-\eqref{ccond} can be rewritten further as a system involving null forms (see \cite{KM, KlainermanTataru, Selberg, MachedonSterbenz, KRT})
 
\begin{align}
\Delta A_0 &= -\Im(\phi \overline{\partial_t\phi} )+|\phi|^2 A_0\tag{MKG-0}\label{e1} \\
\Delta \partial_t A_0 &= -\Im\partial_{j}(\phi \overline{\partial_j\phi} )+\partial_{j}(|\phi|^2 A_j),\tag{MKG-1}\label{e2}\\ 
\Box A_j &= 2R^kD^{-1}Q_{jk} (\Re \phi, \Im \phi) + \mathcal{P}(|\phi|^2 A_j),\tag{MKG-2} \label{e3}\\
\Box \phi &= -iQ_{jk}{(\phi , D^{-1}[R^j A^k -R^k A^j])}+2iA_0 \partial_t \phi + i(\partial_t A_0 ) \phi + A^{\mu}A_{\mu} \phi,\tag{MKG-3} \label{e4}\\
\partial^{j}A_{j}&=0,\tag{MKG-4}\label{e5}
\end{align}
where $R_{k}$ denotes the Riesz transform, $\mathcal P$ is the Leray projection onto the divergence free vector fields, $\mathcal P=\Delta^{-1}d^{\ast}d$, or equivalently
\[
\mathcal P X_{j}=R^{k}(R_{j}X_{k}-R_{k}X_{j}),
\]
and $Q_{\a\b}$ denotes the null form
\[
Q_{\a\b}(u,v)=\partial_{\a}u\partial_{\b}v-\partial_{\b}v\partial_{\a}u.
\]
 
The main result of this paper is contained in the following theorem.

\begin{thm}\label{mainthm}
Let 
$ \frac 14<s'\leq 1,\ \ 
\frac 12<s\leq 1,$ and in addition, let $s'$ satisfy \begin{align}
\max\left(\frac 32 - 2s, \frac s2-\frac 18\right)<s'<4s-\frac 32\label{c2t}.
\end{align}
\begin{figure}
\includegraphics[scale=1]{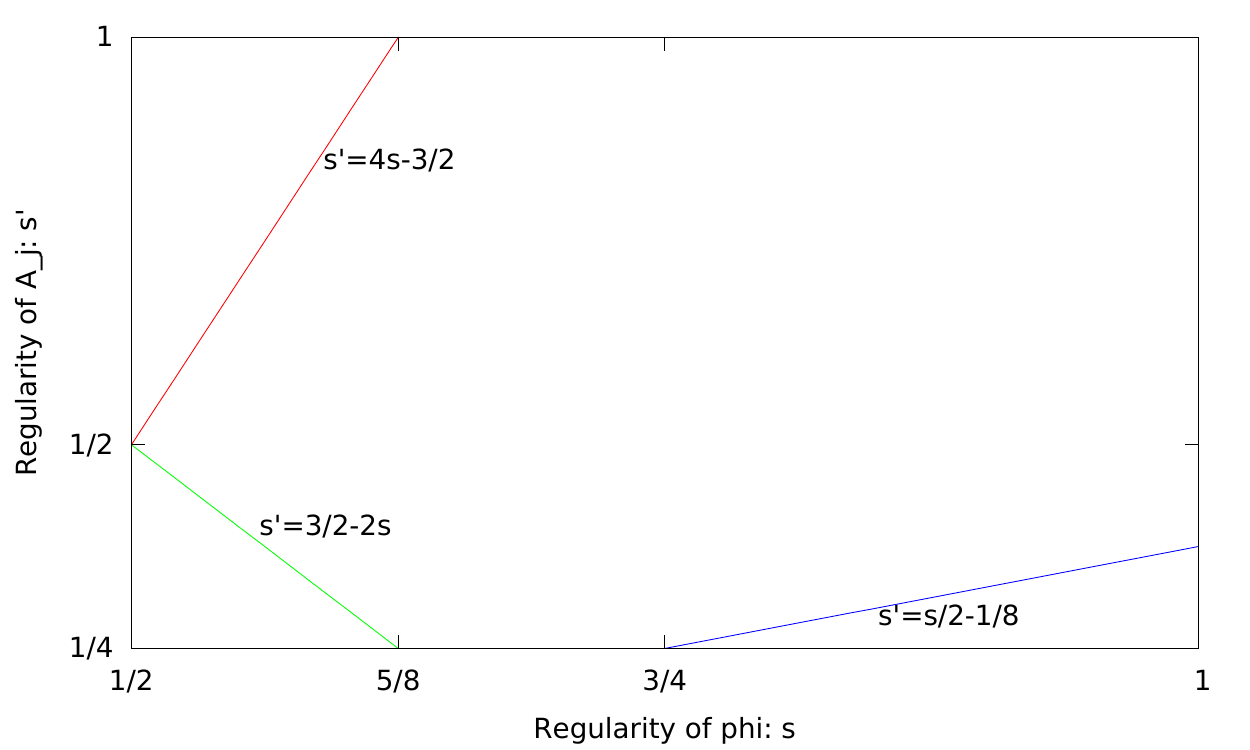}
 \caption{}
\label{figure}
\end{figure}

%\begin {figure}
%\begin{center}
%\input{out}
%\caption{}
%\label{figure}
%\end{center}
%\end {figure}

Consider MKG in the Coulomb gauge given by \eqref{e1}-\eqref{e5} with initial data
\begin{align}
(A_1, A_2,\phi)|_{t=0}&=(a_1,a_2,\phi_0)\in H^{s'}\times H^{s'}\times H^s,\label{id1}\\
(\partial_tA_1, \partial_tA_2,\partial_t \phi)|_{t=0}&=(b_1,b_2,\phi_1)\in H^{s'-1}\times H^{s'-1}\times H^{s-1}\label{id2},\\
\partial^ja_j&=\partial^jb_j=0\label{id3}.
\end{align}

Then the Cauchy problem \eqref{e1}-\eqref{e5}, \eqref{id1}-\eqref{id3} is locally well-posed in the following sense:
\begin{itemize}
\item \textbf{{(Local Existence)}} 
For data given by \eqref{id1}-\eqref{id3} there exist $T>0$ depending continuously on the size of the initial data, and functions 
\begin{align*}
A_0 \in C_b([0,T];\dot H^{\sigma'})&\cap C^1_b([0,T]; \dot H^{\sigma'-1}) ,\ 0<\sigma'< 1+2s,\\
 A_1,A_2\in  C_b([0,T];H^{s'})&\cap C_b^1([0,T];H^{s'-1}) ,\\
   \phi \in C_b([0,T];H^s)&\cap C_b^1([0,T];H^{s-1}),
\end{align*}
which solve (MKG) in the Coulomb gauge on $[0,T] \times \R^2$ in the sense of distributions and such that the initial conditions are satisfied.
\item \textbf{(Uniqueness)} If $T>0$ and $(A,\phi)$ and $(A',\phi')$ are two solutions of (MKG) in the Coulomb gauge on $(0,T)\times \R^2$ belonging to 
\[
(C_b([0,T];\dot H^{\sigma'})\cap C_b^1([0,T];\dot H^{\sigma'-1}) )\times(\mathcal H^{s',\theta}_T)^2\times \mathcal H^{s,\theta}_T
\] 
with the same initial data, then $(A,\phi)=(A',\phi')$ on $(0,T)\times \R^2$.
\item \textbf{(Continuous Dependence on the Initial Data)} For any $(a_1,a_2,\phi_0) \in (H^{s'})^2\times H^{s}$ and $(b_1,b_2,\phi_1) \in (H^{s'-1})^2\times H^{s-1}$ satisfying \eqref{id3} there is a neighborhood $U$ of $(a_1,a_2,\phi_0)\times (b_1,b_2,\phi_1)$ in $ (H^{s'})^2\times H^{s} \times (H^{s'-1})^2\times H^{s-1}$ such that the solution map $$(a_1,a_2, \phi_0) \times (b_1,b_2,\phi_1) \mapsto (A,\phi)$$ is continuous from $U$ into $(C_b([0,T];\dot H^{\sigma'})\cap C_b^1([0,T];\dot H^{\sigma'-1}) )\times (C_b([0,T];H^{s'})\cap C_b^1([0,T];H^{s'-1}))^2\times C_b([0,T]; H^s)\cap C^1_b([0,T]; H^{s-1})$  .  
\end{itemize}
\end{thm}

An immediate corollary (can be also seen from Figure \ref{figure}) is the following

\begin{cor}\label{cor-thm}
 Let  $\frac 12< s\leq 1$.
 Then 2D MKG in the Coulomb gauge is locally well-posed (in the sense stated above) for initial data in $(H^s)^3\times (H^{s-1})^3$.

\end{cor}

\begin{remark}
We do not consider $s', s>1$ since then the initial data is in $L^\infty$ and the estimates are easier.
\end{remark}

\begin{remark}
Figure \ref{figure} shows the region, where we can obtain LWP.  The region is contained between the three lines and bounded below by $\frac 14$.  The region does not include any of the lines except $s,s'=1$.  It allows to take $s'=\frac 14+\epsilon$ for $s\in [\frac 58,\frac 34+\epsilon)$.  After that, the values for $s'$ are bounded below by one of the lines and require $s'>\frac s2-\frac 18$. 
\end{remark}

\begin{remark}
Spaces $\mathcal H^{s,\theta}_T$ are defined in Section \ref{spaces}.
\end{remark}

 \begin{remark}
There is no initial data given for $A_0,$ because $A_0(0)$ can be determined by solving the elliptic equation.  We note though that $T$ depends on the data \eqref{id1}-\eqref{id3} \emph{and on $A(0)$}.  See Section \ref{A0} for more details.  
\end{remark}
The outline of the paper is as follows.  Section \ref{prelim} sets notation, introduces spaces, and estimates used.  In Section \ref{A0} we address the complications that arise in 2D when solving for the elliptic variable $A_{0}$. Section \ref{main} is devoted to the proof of Theorem \ref{mainthm}, which is reduced to establishing appropriate estimates. 

\begin{acknowledgment}
The first author was partially supported by a grant from the Simons Foundation \#246255.
\end{acknowledgment}
\section{Preliminaries}\label{prelim}
First we establish notation, then we introduce function spaces as well as estimates used.  
\subsection{Notation}
We use $a \lesssim b$ to denote $a \leq Cb$ for some positive constant $C$.  Also, $a\simeq b$ means $a\lesssim b$ and $b \lesssim a$.  A point in the $2+1$ dimensional Minkowski space is written as $(t,x)=(x^\alpha)_{0 \leq \alpha \leq 2}.$  We also use $\underline U$ to denote just the spatial part $(U_1, U_2)$ of a vector $(U_0, U_1, U_2)$. Greek indices range from $0$ to $2$, and Roman indices range from $1$ to $2$.  We raise and lower indices with the Minkowski metric, $\diag(-1,1,1)$.  We write $\partial_\alpha=\partial_{x^\alpha}$ and $\partial_t=\partial_0$, and we also use the Einstein notation.  Therefore, $\partial^i \partial_i =\triangle,$ and $\partial^\alpha \partial_\alpha=-\partial^2_t +\triangle=\square$.

 %%%%%%%%%%%%%%%%%%%%%%%%%%%%%%%%%%%%%%%%%%%%%%%%%%%%%%%%%%%%%%%%%%%%%%%%%%%%%%%%%%%%%%%%%%%%%%%%%%%%%%%%%%%%%%%%%%%%%%%%%%%
%%%%%%%%%%%%%%%%%%%				Function Spaces
%%%%%%%%%%%%%%%%%%%%
%%%%%%%%%%%%%%%%%%%%%%%%%%%%%%%%%%%%%%%%%%%%%%%%%%%%%%%%%%%%%%%%%%%%%%%%%%%%%%%%%%%%%%%%%%%%%%%%%%%%%%%%%%%%%%%%%
\subsection{Function Spaces}\label{spaces}
Define following Fourier multiplier operators
\begin{align}
\widehat{\Lambda^\alpha f}(\xi)&=(1 + |\xi|^2)^\frac{\alpha}{2}\hat f(\xi),\\
\widehat{\Lambda^\alpha_+ u}(\tau,\xi)&=(1 + \tau^2+|\xi|^2)^\frac{\alpha}{2}\hat u(\tau,\xi),\\
\widehat{\Lambda^\alpha_- u}(\tau,\xi)&=\left(1 +\frac{ (\tau^2-|\xi|^2)^2}{1+\tau^2+|\xi|^2}\right)^\frac{\alpha}{2}\hat u(\tau,\xi),
\end{align}
where the symbol of $\Lambda^\alpha_-$ is comparable to $(1 + \big||\tau|-|\xi|\big|)^\alpha$.  The corresponding homogeneous operators are denoted by $D^\alpha,D_+^\alpha,D_-^\alpha$ respectively.  
 
 We employ spaces $\h$ and $\mathcal H^{s,\theta}$ with norms given by
\begin{align*}
\|u\|_{H^{s,\theta}}&=\| \Lambda^s\Lambda^\theta_-u \|_{L^2(\R^{2+1})},\\
\|u\|_\H&=\|u\|_{H^{s,\theta}}+\|\partial_tu\|_{H^{s-1,\theta}}.
\end{align*}
An equivalent norm for $\H$ is $\|u\|_{\H}=\| \Lambda^{s-1}\Lambda_+\Lambda^\theta_-u \|_{L^2(\R^{2+1})}$.  If $\theta > \frac{1}{2}$, we have (see for example \cite{SelbergT} )
\begin{align}
H^{s,\theta} &\hookrightarrow  C_b(\R; H^s),\label{embed}\\
\H &\hookrightarrow  C_b(\R; H^s)\cap C_b^1(\R; H^{s-1}).
\end{align}
%This is a crucial fact needed to localize our solutions in time. 
 We denote the restrictions to the time interval $[0,T]$ by
\[
 H^{s,\theta}_T\qand \mathcal H^{s,\theta}_T,
\]
respectively.
%%%%%%%%%%%%%%%%%%%%%%%%%%%%%%%%%%%%%%%%%%%%%%%%%%%%%%%%%%%%%%%%%%%%%%%%%%%%%%%%%%%%%%%%%%%%%%%
%%%%%%%%%%%%%                       Estimates USED
%%%%%%%%%%%%%%%%%%%%%%%%%%%%%%%%%%%%%%%%%%%%%%%%%%%%%%%%%%%%%%%%%%%%%%%%%%%%%%%%%%%%%%%%%%%%%%%
\subsection{{Estimates Used}}  We use two kinds of product estimates.  For Sobolev spaces we have
\be\label{sp}
\norm{uv}_{H^{-s_0}}\lesssim \norm{u}_{H^{s_1}}\norm{v}_{H^{s_2}},
\ee
where $s_0, s_1, s_2$ satisfy $s_0+s_1+s_2\geq 1,$ and $s_0+s_1+s_2\geq \max(s_0, s_1, s_2)$ and at most one of these inequalities is an equality (see for instance \cite{DAFS}).
The $\h$ analog is the following theorem.
\begin{thm}\cite{DAFS}\label{atlas}  
Let $s_0, s_1, s_2, b_0, b_1, b_2 \in \mathbb{R}$, then the following estimate holds for all $u,v \in \mathcal{S}(\R^{2+1})$
$$\|uv\|_{H^{-s_0, -b_0}}\lesssim \|u\|_{H^{s_1, b_1}}\|v\|_{H^{s_2, b_2}}$$
provided that the following conditions are satisfied:
\begin{subequations}
\begin{align}
b_0+b_1+b_2&>\frac{1}{2}\label{a}\\
b_0+b_1&\geq 0\label{b}\\
b_0+b_2&\geq 0\label{c}\\
b_1+b_2&\geq 0\label{d}\\
s_0+s_1+s_2&>\frac{3}{2}-(b_0+b_1+b_2)\label{e}\\
s_0+s_1+s_2&>1-\min(b_0+b_1,b_0+b_2, b_1+b_2)\label{f}\\
s_0+s_1+s_2&>\frac{1}{2}-\min(b_0,b_1,b_2)\label{g}\\
s_0+s_1+s_2&>\frac{3}{4}\label{h}\\
(s_0+b_0)+2s_1+2s_2&>1\label{i}\\
2s_0+(s_1+b_1)+2s_2&>1\label{j}\\
2s_0+2s_1+(s_2+b_2)&>1\label{k}\\
s_1+s_2&\geq \max(0, -b_0)\label{l}\\
s_0+s_2&\geq \max(0, -b_1)\label{m}\\
s_0+s_1&\geq \max(0, -b_2)\label{n}
\end{align} 
\end{subequations}
 \end{thm}
 
%%%%%%%%%%%%%%%%%%%%%%%%%%%%%%%%%%%%%%%%%%%%%%%%%%%%%%%%%%%%%%%%%%%%%%%%%%%%%%%%%%%%%%%%%%%%%%%%
%%%%%%%%%%%%%%%%%%%%%%%%   Elliptic Estimates
%%%%%%%%%%%%%%%%%%%%%%%%%%%%%%%%%%%%%%%%%%%%%%%%%%%%%%%%%%%%%%%%%%%%%%%%%%%%%%%%%%%%%%%%%%%%%%
\section{Elliptic variables $A_{0}$ and $\partial_{t}A_{0}$}\label{A0}
In this section we address the existence, uniqueness and regularity of the elliptic variable $A_{0}$ and its time derivative.   
\subsection{Solving for $A_{0}$}
We discuss three conventional approaches that we were not able to apply to produce estimates on $A_0$.

First, from variational methods, for each $t$ we could obtain existence and uniqueness of $A_0(t)$ in $H^1(\R^2)$ as the minimizer of 
\[
\int_{\R^2}\abs{\nabla A_0}^2
+\abs{D_0\phi}^2 \dd x.
\]
Then $A_0(t)$ would solve 
\begin{align}\label{eA0}
\Delta A_0 =-\Im  (\phi \overline{ D_0 \phi})
\end{align}
as needed.
However, the complications  arise when we would like to arrange $A_0(t)$ into a solution in some space-time norm.  The complications come from the fact that the variational methods do not give us estimates for the $H^1$ norm, and for example, the clever manipulations used in \cite{KRT} to obtain bounds on the \emph{homogeneous} $\dot H^1$ norm only seem to work in 3D or higher (or at a higher regularity in 2D: $s\geq 1$).  Moreover, even in $3D$ the authors were not able to obtain estimates to bound the $L^2$ norm of $A_0$ and had to isolate low frequencies.

Another choice is to resort to the fixed point method, just like in \cite{Czubak10}, and solve the elliptic equation $\eqref{eA0}$ for $A_{0}$, but  in 2D, $\log$ is the fundamental solution of the Laplacian, and so far we have not been able to close the iteration in any Sobolev space.  In \cite{Czubak10}, although the elliptic equation was in 2D, it was essentially using the \emph{derivative} of the fundamental solution of the Laplacian.  Hence, it was more tractable.

We could also try to use the Riesz Representation Theorem (in $H^1$), but then again this would not give us uniform estimates for $A_0(t)$ unless we have $s\geq 1$ (compare with 4D in \cite{Selberg}).  (We allow $s=1$ in this paper, but again we are really interested in $s<1$, so we seek a method that works below $1$.)  
We note however, that we could get an estimate on $\dot H^1$, but in 2D this is not useful unless we include BMO in the estimates.   
%Moreover, the $\dot H^1$ estimate and we would still need to rely on what we know about regularity of the time derivative %$B_0$ as discussed below.

Fortunately, there is another choice.  The equation for $\partial_t A_0$ is better behaved than \eqref{eA0}, and we can solve for $B_{0}=\partial_t A_0$ in $C_b([0,T]; \dot H^\sigma)$ (see Section \ref{B0} below).  Then we can let 
\be\label{A0def}
A_0(t)=\int^t_0 B_0(s)\dd s +a_0,
\ee
where $a_0\in H^1$ is the solution of the variational problem at $t=0$.  Since $\dot H^\sigma \subset \mathcal S'$, \eqref{A0def} defines a tempered distribution.  We need to show that \eqref{A0def} solves the required equation, and that $A_0$ has enough regularity.  In particular, we need $A_0 \in L^\infty$ to handle the estimates for $\phi$ (see Section \ref{ellpiece}), so we will use the equation and bootstrap from the initial $\dot H^\sigma$ estimate (see estimate \eqref{a0sigma}, Lemma \ref{lemmaa0} and Corollary \ref{cora0}).
\begin{remark}
Alternatively, we could argue that we have the existence of the solution in $H^1$ from the variational method or the Riesz Representation Theorem.  Then to obtain estimates on $A_0$, we could show that $B_0$ is the weak time derivative of $A_0$, then use it to show that \eqref{A0def} holds, and then still proceed with \eqref{a0sigma}, Lemma \ref{lemmaa0} and Corollary \ref{cora0}.  
\end{remark}

In \cite{MachedonSterbenz, Selberg} the authors show that if $A_0$ solves \eqref{eA0} and $B_0$ solves \eqref{e2}, then in fact, $B_0=\partial_t A_0$. Here, by definition $\partial_t A_0=B_0$, but it is not immediately obvious that if $A_0$ is defined by \eqref{A0def}, then $A_0$ solves \eqref{eA0}. However, we can show

\lm Let $B_0$ solve \eqref{e2}.  Then $A_0(t)$ given by \eqref{A0def} solves \eqref{eA0} in a sense of tempered distributions  for every $t\in [0,T]$.  
\ml
\begin{proof}
Recall, the current $J=(J_0, \underbar J)$ is given by
\[
J_\alpha(t)=-\Im (\phi(t) \overline {D_\alpha\phi(t)}) ,\quad \alpha=0, 1, 2.
\]
Then \eqref{e2} says
\[
\Delta B_0 =\dv \underbar J,
\]
and we need
\[
\Delta A_0(t)=J_0(t).
\]
Note, it is enough to show the current is conserved, i.e., 
\be\label{cc}
\dv \underbar J= \partial_t J_0,
\ee
because then from \eqref{A0def} we have
 \begin{align*}
 \Delta A_0(t)&=\int^t_0 \Delta B_0(s) \dd s +\Delta a_0\\
 &=\int^t_0 \dv \underbar J(s) \dd s - \Im (\phi (\overline {\partial_t\phi+ia_0\phi})) \\
 &=\int^t_0 \partial_tJ_0(s) \dd s + J_0(0) \\
 &=J_0(t),
 \end{align*}
 as needed.  So we show \eqref{cc}.  To that end, using similar computations as in \cite{MachedonSterbenz, Selberg},  compute
 \be\label{ptj}
 \partial_t J_0=-\Im  (\phi \overline{ \partial_t^2  
\phi}) + (\partial_t \abs{\phi}^2)A_0+\abs{\phi}^2B_0.
 \ee
Then by using \eqref{mkgphi} for $\partial^2_t \phi$, we have
\begin{align*}
\Im  (\phi \overline{ \partial_t^2  
\phi})&=\Im (\phi \overline{2iA^j\partial_j \phi})-\Im (\phi \overline{2i A_0\partial_t \phi})-\Im (\phi \overline{iB_0\phi})+\Im  (\phi \overline{\Delta \phi})\\
&=-\dv \underbar J -\Im (\phi \overline{2i A_0\partial_t \phi})-\Im (\phi \overline{iB_0\phi})\\
&=-\dv \underbar J +(\partial_t \abs{\phi}^2)A_0+ \abs{\phi}^2B_0,
\end{align*}
where to go from the first to the second line, we combined the first and last term using a product rule and the Coulomb condition.  Inserting this into \eqref{ptj} gives \eqref{cc}.
 \end{proof}
 Next we address uniqueness of the solution of \eqref{e1}.
 \lm\label{a0unique} Let $A_0(t)$ be the solution of \eqref{e1}.  Then $A_0(t)$ is unique in $\dot H^\frac 12\cap \dot H^1$.  
%, and $A_0 \in L^\infty([0,T]\times \R^2).$
\ml
\begin{proof}
Let $u, v$ both solve \eqref{e1}.  Then $w=u-v$ solves
\[
-\Delta w+\abs{\phi}^2w=0.
\]
in a sense of tempered distributions.  Because $\mathcal S$ is dense in $\dot H^\frac 12\cap \dot H^1$ this implies
\[
\int_{\R^2}(\abs{\nabla w}^2+\abs{\phi}^2w^2)\dx=0.
\]
so $w=0$ a.e. in $\R^2$.
\end{proof}
 \begin{remark}
 $\dot H^\frac 12$ here is only convenient and not necessary.  We show below that $A_0$ has actually better regularity than just $\dot H^\frac 12\cap \dot H^1$.
   \end{remark}
 \subsection{Solving for $B_{0}$}\label{B0}
Recall \eqref{e2}
\begin{align*} \Delta \partial_t A_0 = -{\partial}_j \Im(\phi \overline{\partial_j\phi} )+{\partial}_j( |\phi|^2 A_j).\end{align*}
So we let
\be\label{B0def}
B_{0}=\frac{\partial_{j}}{\Delta}J_{j}=-2\pi\frac{x_{j}}{\abs{x}^{2}}\ast\left(\Im (\phi\overline{\partial_{j}\phi}) +\abs{\phi}^{2}A_{j}\right),
\ee
and then again, define $A_{0}$ by \eqref{A0def}.  Then \eqref{e2} is satisfied.

\subsection{Estimates for $A_0$ and $\partial_t A_0$}
We start with 
\lm 
%Suppose $\norm{\phi}_C_{b}([0,T]; H^{s}, \norm{A_j}_C_{b}([0,T]; H^{s'}\leq C$, where $C$ is some universal constant
Let $B_{0}(t)$ be given by \eqref{B0def}. Then $B_{0}$ is the unique solution of \eqref{e2} in $C_{b}([0,T]; \dot H^{\sigma})$, for any $\sigma\in (0, 2s-1)$, and
\be\label{B0est}
\norm{B_0}_{C_{b}([0,T];\dot H^{\sigma})}\lesssim \norm{\phi}_{C_{b}([0,T]; H^{s})}^2(1+\norm{\underline A}_{C_{b}([0,T]; H^{s'})}).
\ee
\ml
\begin{proof}  First note that by definition and continuity of the right hand side in $\eqref{B0def}$, $B_0$ solves \eqref{e2} and is continuous in time.  Uniqueness will follow from \eqref{B0est}.

Now, fix $t\in [0,T]$.  Using \eqref{e2} (or \eqref{B0def}), we would like to show 
\begin{align}
\|\phi(t) \overline{{\partial}_j \phi(t)}\|_{\dot H^{\sigma-1}}&\lesssim \norm{\phi(t)}_{H^{s}}\norm{\partial_{j}\phi(t)}_{H^{s-1}},\label{b01}\\
 \|\abs{\phi(t)}^2 A_j(t)\|_{\dot{H}^{\sigma-1}}&\lesssim\|\phi(t)\|^2_{H^s} \|A_j(t)\|_{H^{s'}}\label{b02}.
\end{align}

For \eqref{b01}, we use duality and show   \begin{align}\label{diditalready}
\|uv\|_{H^{1-s}}\lesssim \|u\|_{H^s} \|v\|_{\dot{H}^{1-\sigma}}.\end{align}
First 
 \[
 \begin{split}
 \|uv\|_{H^{1-s}}&\lesssim \|uv\|_{L^2} + \|D^{1-s}(uv)\|_{L^2}\label{3terms}\\
 &\lesssim \|uv\|_{L^2} +\|(D^{1-s}u)v\|_{L^2}+ \|uD^{1-s}v\|_{L^2}\\
 &=I+II+III.
 \end{split}
 \]
By H\"older's with $ \frac{1}{2}= {\frac{1-\sigma}{2}}+(\frac{1}{2}-{\frac{1-\sigma}{2}})$ we have
 \begin{align*}
I\lesssim \|u\|_{(\frac 12-\frac \sigma 2)^{-1}}\|v\|_{(\frac{1}{2}-{\frac{1-\sigma}{2}})^{-1}}\lesssim  \|u\|_{\dot{H}^{\sigma}}\norm{v}_{\dot H^{1-\sigma}}\lesssim\|u\|_{H^s}\norm{v}_{\dot H^{1-\sigma}},
 \end{align*}
as long as 
\be\label{sigma1}
\sigma \leq s\quad\mbox{ and}\quad 0<\sigma<1.
\ee

Next, by the same application of H\"older's
 \begin{align*}
II\lesssim \|D^{1-s}u\|_{(\frac 12-\frac \sigma 2)^{-1}}\|v\|_{(\frac{1}{2}-{\frac{1-\sigma}{2}})^{-1}}\lesssim  \|D^{1-s}u\|_{\dot{H}^{\sigma}}\norm{v}_{\dot H^{1-\sigma}}\lesssim\|u\|_{H^s}\norm{v}_{\dot H^{1-\sigma}},
 \end{align*}
 provided
 \be\label{sigma2}
\sigma\leq 2s-1\quad\mbox{ and}\quad 0<\sigma<1.
\ee
Similarly, by H\"older's with  $\frac{1}{2}=\frac{s-\sigma}{2}+(\frac{1}{2}-\frac{s-\sigma}{2})=\frac 12-\frac{1-s+\sigma}{2}+(\frac{1}{2}-\frac{s-\sigma}{2})$

 \begin{align*}III&\lesssim \|u\|_{(\frac 12-\frac{1-s+\sigma}{2})^{-1}}\|D^{1-s}v\|_{(\frac{1}{2}-\frac{s-\sigma}{2})^{-1}}\\
 &\lesssim \|u\|_{\dot H^{1-s+\sigma} }\|D^{1-s+s-\sigma}v\|_{2}\\
 &\lesssim \|u\|_{H^{s}}\|v\|_{\dot H^{1-\sigma}}
  \end{align*}
if
 \be\label{sigma3}
\sigma\leq 2s-1\quad\mbox{ and}\quad 0<\sigma<s.
\ee

This completes the proof of \eqref{b01}.  Now for \eqref{b02}, we would like to show 
\begin{align*}
\|uvw\|_{\dot{H}^{\sigma-1}}\lesssim\|u\|_{H^s}\|v\|_{H^s}\|w\|_{H^{s'}}.\end{align*}

By Sobolev embedding and H\"older's we have  
\begin{align*}\|uvw\|_{\dot{H}^{\sigma-1}}&\lesssim\|uvw\|_{\frac{2}{2-\sigma}}\\
 &\lesssim\|uv\|_{{\frac{2}{1-\sigma}}}\|w\|_2\\
  &\lesssim\|uv\|_{ H^\sigma}\|w\|_{H^{s'}}\\
    &\lesssim\|u\|_{ H^s}\|v\|_{ H^s}\|w\|_{H^{s'}},
 \end{align*}
 where to go to the last line we use \eqref{sp} provided
\be\label{sigma4}
\sigma\leq 2s-1\quad\mbox{ and}\quad \sigma<s.
\ee
Collecting \eqref{sigma1}-\eqref{sigma4} and using $s\leq 1$, we obtain $\sigma\in (0, 2s-1).$

\end{proof}
  
Next, from \eqref{A0def} and \eqref{B0est}, we immediately get $A_0 \in C_b([0,T]; \dot H^\sigma)$ and
\begin{align}\label{a0sigma}
\norm{A_0}_{C_{b}([0,T];  \dot H^\sigma)}&\lesssim T\norm{\phi}_{C_{b}([0,T]; H^{s})}^2(1+ \norm{\underline A}_{C_{b}([0,T]; H^{s'})})+\norm{a_0}_{\dot H^\sigma},
 \end{align}
where $0<\sigma< 2s-1$.
Note, from the variational method we have $a_0 \in H^1$, so $\norm{a_0}_{\dot H^\sigma}$ is finite; we just do not have the estimates to control it in terms of the data for $\phi$ and $\underline A$.   Also, because $\norm{a_0}_{\dot H^\sigma}$ appears in \eqref{a0sigma}, it will appear in \eqref{elliptic1} and \eqref{elliptic3}, and hence $T$ will also depend on $\norm{a_0}_{\dot H^\sigma}$.
\lm\label{lemmaa0} Let $\frac 12<s\leq 1$, and $A_0 \in C_b([0,T]; \dot H^\sigma),\ 0<\sigma < 2s-1$.  Then $A_0 \in C_b([0,T]; \dot H^a)$, where $1<a< 2s$, and
\[
\norm{A_0}_{C_b([0,T]; \dot H^a)}\lesssim \norm{\phi}_{C_{b}([0,T]; H^{s})}\norm{\phi_t}_{C_{b}([0,T]; H^{s-1})}+ \norm{\phi}_{C_{b}([0,T]; H^{s})}^2\norm{A_0}_{C_{b}([0,T]; \dot H^{\sigma})},
\]
where $0<\sigma<2s-1$.
\ml
\begin{proof}  We have
\[
\norm{D^a A_0(t)}_2=\norm{D^{a-2}\Delta A_0(t)}_2=\norm{D^{a-2}J_0}_2.
\]
So we need to estimate $\norm{\phi \phi_t}_{\dot H^{a-2}}$ and $\norm{A_0 \abs{\phi}^2}_{\dot H^{a-2}}$.
For the first estimate we need
\begin{align*}
H^s \cdot H^{s-1}\hookrightarrow \dot H^{a-2},
\end{align*}
which is equivalent by duality to
 \begin{align*}
\norm{uv}_{H^{1-s}}\lesssim \norm{u}_{H^s}\norm{v}_{\dot H^{2-a}},
\end{align*}
but from the assumptions on $a$ and $\sigma$, this is exactly the estimate \eqref{diditalready}.

So we bound the cubic term.  Using Sobolev with $\frac12=\frac 1p-\frac{2-a}{2}$
\[
\norm{A_0 \abs{\phi}^2}_{\dot H^{a-2}}\lesssim \norm{A_0\abs{\phi}^2}_p.
\]
We will be done by H\"older and Sobolev, if we can write $\frac 1p=(\frac 12-\frac\sigma 2)+2(\frac 12-\frac \alpha2 )$ for some $0<\alpha\leq s$, and where $\sigma$ is the number of the derivatives we have on $A_0$ using \eqref{a0sigma}.  But 
\[
\frac 1p=\frac 12+\frac{2-a}{2}=\left(\frac 12-\frac\sigma2\right)+\frac{2-a+\sigma}{2}=\left(\frac 12-\frac\sigma2\right)+2\left(\frac 12-\frac{\tfrac a2-\tfrac \sigma 2}{2}\right),
\]
and $\frac a2-\frac \sigma 2\leq s$ as needed.
\end{proof}
By interpolation using \eqref{a0sigma} and Lemma \ref{lemmaa0} we have
\begin{cor}\label{cora0} $A_0 \in C_b([0,T]; \dot H^a), \ 0<a< 2s$.
%, and $A_0 \in L^\infty([0,T]\times \R^2).$
\end{cor}  
\begin{remark} We note the difference in one derivative on the estimates for $A_0$ and $B_0$.  Since $B_0$ is the time derivative of $A_0$, this difference on the spatial estimates is quite natural.
\end{remark}

%%%%%%%%%%%%%%%%%%%%%%%%%%%%%%%%%%%%%%%%%%%%%%%%%%%%%%%%%%%%%%%%%%%%%%%%%%%%%%%%%%%%%%%%%%%%%%%%%%%%%%%%%%%%%%%%
%%%%%%%%%%%%%%            PROOF OF THE MAIN THEOREM
%%%%%%%%%%%%%%%%%%%%%%%%%%%%%%%%%%%%%%%%%%%%%%%%%%%%%%%%%%%%%%%%%%%%%%%%%%%%%%%%%%%%%%%%%%%%%%%%%%%%%%%%%%%%%%%%%
\section{Proof of Theorem \ref{mainthm}}\label{main}
As is now well-known (see for example \cite{SelbergT, SelbergEstimates}), to show Theorem \ref{mainthm} it is enough to estimate the nonlinearities in the appropriate spaces.  
 %%%%%%%%%%%%%%%%%%%%%%%%%%%%%%%%%%%%%%%%%%%%%%%%%%%%%%%%%%%%%%%%%%%%%%%%%%%%%%%%%%%%%%%%%%%%%%%%%%%%%%%%%%%%%%%%%%%%%%%%%%%%%%%%%%%%%%%%%%%
%%%%%%%%%%%%%%%%%%%%%%%%%%%%%                    SET UP OF THE ITERATION
%%%%%%%%%%%%%%%%%%%%%%%%%%%%%%%%%%%%%%%%%%%%%%%%%%%%%%%%%%%%%%%%%%%%%%%%%%%%%%%%%%%%%%%%%%%%%%%%%%%%%%%%%%%%%%%%%%%%%%%%%%%%%%%%%%%%%%%%%%%
 
\begin{subsection}{Estimates needed}
The estimates for the elliptic equations are discussed in Section \ref{A0}.   For the wave equations we need to estimate
\begin{align*}
\norm{\Lambda^{-1}_+\Lambda^{-1+\epsilon}_-\mathcal \Box A_j^{(m)}}_{\Ha}&=\norm{\Box A_j^{(m)}}_{H^{s'-1,\theta_1-1+\eps}} \\
\intertext{and}
\|\Lambda^{-1}_+\Lambda^{-1+\epsilon}_-\mathcal \Box \phi^{(m)}\|_{\h}&=\norm{\Box \phi^{(m)}}_{H^{s-1,\theta_0-1+\eps}} .
 \end{align*}

 Since the Riesz transforms are clearly bounded on $L^2$, and the Leray projection $\mathcal P$ is defined in terms of Riesz transforms, we ignore them in the estimates needed.  So it is enough to prove the following
\begin{align}
 \|D^{-1}Q_{jk} (\Re \phi, \Im \phi)\|_{H^{s'-1,\theta_1-1+\epsilon}}&\lesssim \|\phi\|_{H^{s,\theta_0}}^2,\label{nullA}\\
\||\phi|^2 A_j\|_{H^{s'-1,\theta_1-1+\epsilon}}&\lesssim  \|\phi\|_{H^{s,\theta_0}}^2\|A_j\|_{H^{s',\theta_1}}, \label{cubicA}\\
\|Q_{jk}{(\phi , D^{-1}A_j )}\|_{H^{s-1,\theta_0-1+\epsilon}} &\lesssim \|\phi\|_{H^{s,\theta_0}}\|A_j\|_{H^{s',\theta_1}} ,\label{nullphi}\\
 \|A_0 \partial_t \phi\|_{H^{s-1,\theta_0-1+\epsilon}} &\lesssim \|A_0\|_{X_0} \|\partial_t\phi\|_{H^{s-1,\theta_0}},\label{elliptic1}\\
 \|\partial_tA_0 \phi\|_{H^{s-1,\theta_0-1+\epsilon}} &\lesssim \|\partial_tA_0\|_{C_b([0,T]; \dot H^\sigma)} \|\phi\|_{H^{s,\theta_0}},  \ 0<\sigma<2s-1,\label{elliptic2}\\
\|A^2_j \phi\|_{H^{s-1,\theta_0-1+\epsilon}} &\lesssim \|\phi\|_{H^{s,\theta_0}}\|A_j\|_{H^{s',\theta_1}}^2,\label{cubicphi}\\
\|A^2_0 \phi\|_{H^{s-1,\theta_0-1+\epsilon}} &\lesssim \|\phi\|_{H^{s,\theta_0}}\|A_0\|_{X_0}^2\label{elliptic3},
\end{align}
where we let
\begin{align*}
X_0&=C_b([0,T]; L^\infty\cap \dot H^a), \ 1<a<2s.
\end{align*}
Given $(s,s')$ such that
\begin{align}
\frac 14<s'\leq 1,\quad
\frac 12<s\leq 1,\label{c1}
\end{align}
and in addition
\begin{align}
\max\left(\frac 32 - 2s, \frac s2-\frac 18\right)<s'<4s-\frac 32\label{c2},
\end{align}
choose $\theta_0, \theta_1$ such that
\begin{align}
%\max \left(\frac 12,  \frac 18+\frac {\theta_1}{ 2},  \theta_1-\frac 14 \right)<&\theta_0<1,\\
\frac 12<&\theta_0<\frac 34,\\
\max\left(\frac 12, 1-s'\right)  <&\theta_1<\min \left(\frac 34, 4s-s'-1, 2s -\frac 12 \right)\label{theta_1}.
\end{align}
The restrictions on $(s, s')$ allow us to find $\theta_1$ satisfying the required conditions.  For convenience of the reader, we add that at some point it is needed that $s'<2s-\frac 14$, but that is guaranteed by the current upper bound on $s'$.

%%%%%%%%%%%%%%%%%%%%%%%%%%%%%%%%%%%%%%%%%%%%%%%%%%%%%%%%%%%%%%%%%%%%%%%%%%%%%%%%%%%%%%%%%%%%%%%%%%%%%%%%%%%%%%%%%%%
%%%%%%%%%%%%%%%%%%%%                            NULL FORMS                                      %%%%%%%%%%%%%%%%%%%
%%%%%%%%%%%%%%%%%%%%%%%%%%%%%%%%%%%%%%%%%%%%%%%%%%%%%%%%%%%%%%%%%%%%%%%%%%%%%%%%%%%%%%%%%%%%%%%%%%%%%%%%%%%%%%%%%%%
\subsubsection{Null Forms--Proof of Estimate \eqref{nullA}}\label{nullformqij}
\eqref{nullA} will follow from
\be\label{na1}
 \|D^{-1}Q_{jk} (u,v)\|_{H^{s'-1,\theta_1-1+\epsilon}}\lesssim \|u\|_{H^{s,\theta_0}}\|v\|_{H^{s,\theta_0}}.
\ee
Start by  recalling \cite{KM3}
\be\label{qij}
Q_{jk} (u,v)\precsim D^{\frac{1}{2}}{D_{-}}^{\frac{1}{2}}(D^{\frac{1}{2}}uD^{\frac{1}{2}}v)
+ D^{\frac{1}{2}}({D_{-}}^{\frac{1}{2}}D^{\frac{1}{2}}uD^{\frac{1}{2}}v)+D^{\frac{1}{2}}(D^{\frac{1}{2}}u{D_{-}}^{\frac{1}{2}}D^{\frac{1}{2}}v),
\ee
where $w_{1} \precsim w_{2}$ means $|\hat w_{1}|\leq C \hat w_{2}$ for some $C>0$. 
 
Next, following \cite{SelbergK}, we estimate $D^{-1}$ by
\be\label{lohi}
D^{-1}\precsim \Lambda^{-1}|_{\{\abs{\xi}\geq 1\}}+\Lambda^{-M}D^{-1}|_{\{\abs{\xi}<1\}},
\ee
where $M>0$ can be taken as large as we wish.
Then for high frequencies using \eqref{qij} 
since $u$ and $v$ have the same regularity, \eqref{na1} follows from showing
 
\begin{align}
 \|uv\|_{H^{s'-\frac{3}{2},\theta_{1}-\frac{1}{2}+\epsilon}}&\lesssim \|u\|_{H^{s-\frac{1}{2}, \theta_0}}\|v\|_{H^{s-\frac{1}{2}, \theta_0}},\\ 
 \|uv\|_{H^{s'-\frac{3}{2},\theta_{1}-1+\epsilon}}&\lesssim \|u\|_{H^{s-\frac{1}{2}, \theta_0-\frac{1}{2}}}\|v\|_{H^{s-\frac{1}{2}, \theta_0}}.
\end{align}

Then given \eqref{c1}-\eqref{theta_1} the estimates hold by Theorem \ref{atlas}.
 
Now for low frequencies instead of \eqref{qij} we can use a simpler estimate \cite[p. 272]{SelbergK}
\[
Q_{ij}(u,v)\precsim D(D^{\frac 12}u D^{\frac 12}v).
\]
This reduces \eqref{na1} to
 \be\label{likeb4}
 \|uv\|_{H^{-M,\theta_1-1+\epsilon}}\leq \norm{uv}_{H^{-M,0}}\lesssim\|u\|_{L^2_tH^{s-\frac 12}} \|v\|_{L^\infty_tH^{s-\frac 12}}\lesssim \|u\|_{H^{s-\frac 12,\theta_0}} \|v\|_{H^{s-\frac 12,\theta_0}} ,
 \ee
where the first inequality holds because $\theta_1-1+\epsilon\leq 0$ and the third one follows from \eqref{embed} and the trivial embedding $\norm{u}_{H^{0,0}}\lesssim \norm{u}_{H^{0,\alpha}}$ for any $\alpha\geq 0$.  Finally, the second inequality follows from a spatial estimate
\[
\norm{uv}_{H^{-M}}\lesssim\norm{u}_{H^{s-\frac 12}}\norm{v}_{H^{s-\frac 12}}.
\]
which in turn holds by \eqref{sp} if $M$ is large enough and $s>\frac 12$.
%%%%%%%%%%%%%%%%%%%%%%%%%%%%%%%%%%%%%%%%%%%%%%%%%%
%%%%%%%% cubic for A
%%%%%%%%%%%%%%%%%%%%%%%%%%%%%%%%%%%%%%%
\subsubsection{Cubic term: Proof of Estimate \eqref{cubicA}} 
 Again for convenience we record the estimate which implies \eqref{cubicA}
\be\label{cubicAa}
\|uvw\|_{H^{s'-1,\theta_1-1+\epsilon}}\lesssim\|u\|_{H^{s',\theta_1}}\|v\|_{H^{s,\theta_0}}\|w\|_{H^{s,\theta_0}}.
\ee
 The estimate follows easily from two applications of  \eqref{sp} since $s>\frac 12$ ($0<\delta<<1$ in the third term appears to cover the case $s'=1$)
\be\label{cubicAb}
\begin{split}
\|uvw\|_{H^{s'-1,\theta_1-1+\epsilon}}\leq\|uvw\|_{H^{s'-1,0}}&\lesssim \norm{u}_{L^2_tH^{s'}}\norm{uv}_{L^\infty_tH^{\delta}}\\
&\lesssim \norm{u}_{H^{s',\theta_1}}\norm{w}_{H^{s,\theta_0}}\norm{v}_{H^{s,\theta_0}}.
\end{split}
\ee
(Note this holds for any choice of $\frac 12<\theta_0,\theta_1<1$.) 
\begin{remark}\label{rq}

 In fact, if we use Theorem \ref{atlas} we can have both $s$ and $s'$ close to $\frac 14$.  We show this below when we establish \eqref{cubicphi}, which is the estimate \eqref{cubicAa} with the roles of $s, \theta_0$ and $s',\theta_1$ reversed.

\end{remark}

%%%%%%%%%%%%%%%%%%%%%%%%%%%%%%%%%%%%%%%%%%%%%%%%%%
%%%%%%%% null form for phi
%%%%%%%%%%%%%%%%%%%%%%%%%%%%%%%%%%%%%%%

\subsubsection{Null Forms--Proof of Estimate \eqref{nullphi}}\label{nullformq}  We show
\be\label{nullphia}
\|Q_{jk}{(u , D^{-1}v )}\|_{H^{s-1,\theta_0-1+\epsilon}} \lesssim \|u\|_{H^{s,\theta_0}}\|v\|_{H^{s',\theta_1}} 
\ee
We use \eqref{lohi} again.  For high frequencies, \eqref{qij} reduces \eqref{nullphia} to proving that the following three estimates hold
\begin{align}
\|uv\|_{H^{s-\frac{1}{2},\theta_0-\frac{1}{2}+\epsilon}}&\lesssim\|u\|_{H^{s-\frac{1}{2},\theta_0}}\|v\|_{H^{s'+\frac{1}{2},\theta_1}},\label{nullphi1}\\
\|uv\|_{H^{s-\frac{1}{2},\theta_0-1+\epsilon}}&\lesssim\|u\|_{H^{s-\frac{1}{2},\theta_0-\frac{1}{2}}}\|v\|_{H^{s'+\frac{1}{2},\theta_1}},\\
\|uv\|_{H^{s-\frac{1}{2},\theta_0-1+\epsilon}}&\lesssim\|u\|_{H^{s-\frac{1}{2},\theta_0}}\|v\|_{H^{s'+\frac{1}{2},\theta_1-\frac{1}{2}}}\label{nullphi3}.
\end{align}

Then in view \eqref{c1}-\eqref{theta_1}, the estimates hold by Theorem \ref{atlas}.

For low frequencies instead of \eqref{qij} we use \cite[p. 272]{SelbergK}
\[
Q_{ij}(u,v)\precsim D^\frac 12(D^\frac 12 u Dv),
\]
and reduce \eqref{nullphia}, just like in \eqref{likeb4}, to showing
\[
\norm{uv}_{H^{s-\frac 12}}\lesssim \norm{u}_{H^{s-\frac 12}}\norm{v}_{H^M},
\]
which holds by \eqref{sp} if $M$ is large enough.
%%%%%%%%%%%%%%%%%%%%%%%%%%%%%%%%%%%%%%%%%%%%%%%%%%%%%%%%%%%%%%%%%%%%%%%%%%%%%%%%%%%%%%%%%%%%%%%%%%%%%%%%%%%%%%%%%
%%%%%%%%%%%%%%            PROOF FOR THE ELLIPTIC PIECE
%%%%%%%%%%%%%%%%%%%%%%%%%%%%%%%%%%%%%%%%%%%%%%%%%%%%%%%%%%%%%%%%%%%%%%%%%%%%%%%%%%%%%%%%%%%%%%%%%%%%%%%%%%%%%%%%%
\subsubsection{Elliptic Piece: Proof of Estimate (\ref{elliptic1})}\label{ellpiece}
Recall we wish to show
\be
\|A_0\partial_t\phi\|_{H^{s-1,\theta-1+\epsilon}}\lesssim \norm{A_{0}}_{X_0} \|\partial_t\phi\|_{H^{s-1,\theta_0}},\label{M21}
\ee
where
\begin{align*}
X_0= C_b([0,T]; L^\infty \cap \dot H^a), \ 1<a<2s.
%\\ X_1&=C_b([0,T]; \dot H^\sigma), \ 0<\sigma<2s-1.
\end{align*}

If we could estimate $A_0 \in H^{\sigma,0}$ (just like for example authors did in \cite{KRT}), the left hand side of \eqref{M21} could be bounded using Theorem \ref{atlas} as long as
$\sigma> 0$.  In 2D we have to work a little harder.

Using $\theta-1+\epsilon < 0$ first reduce \eqref{M21} to
\be
\|\ip{D}^{s-1}(A_{0}\partial_t\phi) \|_{L^2(\R^{2+1})}\lesssim \norm{A_{0}}_{X_0} \|\partial_t\phi\|_{H^{s-1,\theta_0}},
\ee
which by H\"older in time can follow from 
\be
\|\ip{D}^{s-1}(uv) \|_{2}\lesssim \norm{u}_{X_{0,x}}\|v\|_{H^{s-1}}\label{m21a},
\ee
where we denote
\[
X_{0,x}=L^\infty \cap \dot H^a, \ 1<a< 2s.
\]
Observe, using Fourier transform we can show $\dot H^{1+\tilde\epsilon}\cap \dot H^{1-\tilde\epsilon}\hookrightarrow L^\infty(\R^2)$ for $\tilde\epsilon >0$, so by Corollary \ref{cora0}, $A_0(t) \in L^\infty$.

Next, by duality \eqref{m21a} is equivalent to
\be
\|uv \|_{H^{1-s}}\lesssim \norm{u}_{X_{0,x}}\|v\|_{H^{1-s}}\label{m21b}.
\ee

 We can suppose $\frac 12<s<1$ since if $s=1$, \eqref{m21b} follows by H\"older.  In this case we restrict $X_{0,x}$ to be
\[
X_{0,x}=L^\infty \cap \dot H^a, \ 1<a< \min(2s, 2-s).
\]
We estimate
\[
\|uv \|_{H^{1-s}}\lesssim \norm{uv}_2+\norm{(D^{1-s}u)v}_2+\norm{u(D^{1-s}v)}_2,
\]
which we can bound as follows
\begin{align*}
 \norm{uv}_2&\leq \norm{u}_\infty\norm{v}_2,\\
 \norm{u(D^{1-s}v)}_2&\leq\norm{u}_\infty\norm{v}_{\dot H^{1-s}}.
 \end{align*}
 Finally
 \begin{align*}
 \norm{(D^{1-s}u)v}_2&\leq \norm{u}_{\dot W^{1-s,p}}\norm{v}_{(\frac 12-\frac 1p)^{-1}}\lesssim \norm{u}_{\dot W^{1-s,p}}\norm{v}_{H^{1-s}},
 \end{align*}
 where the last estimate follows by Sobolev embedding with $\tfrac 12-\frac 1p=\frac 12-\frac{2-a-s}{2}$, and $2-a-s\leq 1- s$ if $p$ is chosen so that
 \be\label{pacond}
p>2, \quad\frac 1p=\frac{2-a-s}{2},\quad 1<  a<\min( 2s,2-s).
\ee
Finally, another application of the Sobolev embedding completes the proof of \eqref{m21b} since
\[
 \norm{u}_{\dot W^{1-s,p}}\lesssim \norm{u}_{\dot H^{a}},\quad \frac 1p=\frac 12-\frac{a-1+s}{2}.
\]

\subsubsection{Elliptic Piece: Proof of Estimate \eqref{elliptic2}}\label{ellpiecet}
Next we need
\be
\|\partial_tA_0\phi\|_{H^{s-1,\theta-1+\epsilon}} \lesssim\|\partial_tA_0\|_{C_b\dot H_{x}^{\sigma}} \|\phi\|_{H^{s,\theta_0}}.\label{M22}
\ee
 It is easy to to see that
\begin{align*}\|\phi\partial_tA_0\|_{H^{s-1,\theta_0-1-\epsilon}}\lesssim\|\phi\partial_tA_0\|_{H^{s-1,0}}=\|\phi\partial_tA_0\|_{L^2_tH^{s-1}}.
\end{align*}
By H\"older in time it is enough to show 
\begin{align*}
\|\phi \partial_tA_0\|_{H^{s-1}}\lesssim \|\phi\|_{H^{s}}\|A_0\|_{\dot{H}^{\sigma}}.
\end{align*}
Then applying Sobolev embedding $\norm{u}_{2}\lesssim \norm{\Lambda^{1-s}u}_{p}$, $\frac 12=\frac 1p-\frac{1-s}{2}$, we get
 \begin{align*}\|\phi\partial_tA_0\|_{H^{s-1}}\lesssim\|\phi \partial_tA_0\|_{L^{\frac{2}{2-s}}}.\end{align*}
By H\"older's inequality we have
\begin{align*}
\|\phi \partial_tA_0\|_{L^{\frac{2}{2-s}}}\lesssim\|\phi\|_{L^\frac{2}{1-s+\sigma}}\|\partial_t A_0\|_{L^{\frac{2}{1-\sigma}}} ,
\end{align*}
and another application of Sobolev embedding gives us
\begin{align*}
\|\phi\|_{L^\frac{2}{1-s+\sigma}}\|\partial_t A_0\|_{L^\frac{2}{1-\sigma}}\lesssim \|\phi\|_{H^{s-\sigma}}\|\partial_t A_0\|_{\dot{H}^{\sigma}}\lesssim \|\phi\|_{H^{s}}\|\partial_t A_0\|_{\dot{H}^{\sigma}}  .
 \end{align*}
\subsubsection{Cubic piece: Proof of Estimate \eqref{cubicphi}}  This is equivalent to showing
\be\label{cubicphia}
\|uvw\|_{H^{s-1,\theta_0-1+\epsilon}}\lesssim \|u\|_{H^{s,\theta_0}}\|v\|_{H^{s',\theta_1}}\|w\|_{H^{s',\theta_1}}\ee
As mentioned in Remark \ref{rq} \eqref{cubicphia} is \eqref{cubicAa} with the roles of $s, \theta_0$ and $s',\theta_1$ switched.  Hence if $s'>\frac 12$, the estimate follows just like in \eqref{cubicAb}.  So we can suppose $\frac 14<s'\leq \frac 12$.  Here we show the estimate for $\frac 14< s\leq 1$. (Assuming $s>\frac 12$ would also not simplify the presentation since we want $s'$ close to $\frac 14$.)  Then we have by Theorem \ref{atlas}
\be\label{cubicphib}
\|uvw\|_{H^{s-1,\theta_0-1+\epsilon}}\lesssim \|u\|_{H^{s,\theta_0}}\|vw\|
_{H^{2s'-\frac{3}{4}-\epsilon,0}},
\ee
provided $s<2s'+\frac{1}{4}$.
Another iteration of Theorem \ref{atlas} gives
$$\|vw\|_{H^{2s'-\frac{3}{4}-\epsilon,0}}\lesssim \|v\|_{H^{s',\theta_1}}\|w\|_{H^{s',\theta_1}},$$
and \eqref{cubicphia} follows as needed.

\subsubsection{Elliptic Piece: Proof of Estimate \eqref{elliptic3}}\label{ellpiece3}
This is clear since
\be
\|A^2_0\phi\|_{H^{s-1,\theta_0-1+\epsilon}} \leq \|A^2_0\phi\|_{L^2_{t,x}} \leq \norm{A_{0}}^2_\infty \|\phi\|_{H^{s,\theta_0}}.
\ee

\end{subsection}

%%%%%%%%%%%%%%%%%%%%%%%%%%%%%%%%%%%%%%%%%%%%%%%%%%%%%%%%%%% 
\bibliography{mkgbib}
\bibliographystyle{plain}
\vspace{.125in}
 
\end{document}